# Expansions of the solutions to the confluent Heun equation in terms of the Kummer confluent hypergeometric functions

T.A. Ishkhanyan[1,2] and A.M. Ishkhanyan[1]

[1]Institute for Physical Research, NAS of Armenia, 0203 Ashtarak, Armenia
[2]Moscow Institute of Physics and Technology, Dolgoprudny, Moscow Region, 141700 Russia

**Abstract.** We examine the series expansions of the solutions of the confluent Heun equation in terms of three different sets of the Kummer confluent hypergeometric functions. The coefficients of the expansions in general obey three-term recurrence relations defining double-sided infinite series; however, four-term and two-term relations are also possible in particular cases. The conditions for left- and/or right-side termination of the derived series are discussed.

## 1. Introduction

Expansions of the solutions of the confluent Heun equation [1-3] in terms of mathematical functions other than powers have been discussed by many authors (see, e.g., [3-11]). The Gauss hypergeometric functions [3-5], Kummer and Tricomi confluent hypergeometric functions [6-8], Coulomb wave functions [8,9], Bessel and Hankel functions [10], incomplete Beta functions [11], and other standard special functions have been applied as expansion functions. Using the properties of the derivatives of the solutions of the Heun equation [12-16], it is possible to construct expansions in terms of higher transcendental functions, e.g., the Goursat generalized hypergeometric functions [16] and the Appell generalized hypergeometric functions of two variables of the first kind [17]. Here we discuss several expansions in terms of the Kummer confluent hypergeometric functions (confluent hypergeometric functions of the first kind) starting from the differential equation and the recurrence relations the latter functions obey for chosen forms of the dependence on the summation variable. In general, the coefficients of the expansions obey three-term recurrence relations; however, for one of the discussed forms of involved confluent hypergeometric functions a four-term recurrence relation is also possible. Besides, for a specific choice of the involved parameters a different two-term recurrence relation is obtained. As a result, the expansion coefficients in this case are explicitly calculated in terms of the Gamma functions. Since the forms of used Kummer confluent hypergeometric functions differ from those applied in previous discussions the conditions for termination of the presented expansions refer to different choices of the involved parameters.

The confluent Heun equation is a second order linear differential equation having regular singularities at $z=0$ and $1$, and an irregular singularity of rank 1 at $z=\infty$. We adopt here the following form of this equation [1]:

$$u'' + \left(\frac{\gamma}{z} + \frac{\delta}{z-1} + \varepsilon\right)u' + \frac{\alpha z - q}{z(z-1)}u = 0, \qquad (1)$$

which slightly differs from that adopted in [3] since the parameters $\varepsilon$ and $\alpha$ are here assumed to be independent. This is a useful convention for practical applications since in this form the equation includes the Whittaker-Ince limit [18] of the confluent Heun equation as a particular case achieved by the simple choice $\varepsilon=0$. Note, however, that the expansions presented below do not apply to this limit. The corresponding confluent hypergeometric expansions in terms of the Kummer and Tricomi functions for this case are discussed in [19].

We search for expansions of the solutions of Eq. (1) in the form

$$u = \sum_n a_n u_n, \quad u_n = {}_1F_1(\alpha_n; \gamma_n; s_0 z), \qquad (2)$$

where ${}_1F_1(\alpha_n;\gamma_n;s_0 z)$ is the Kummer confluent hypergeometric function (confluent hypergeometric function of the first kind). Functions $u_n$ obey the confluent hypergeometric differential equation

$$u_n'' + \left(\frac{\gamma_n}{z} - s_0\right)u_n' - \frac{\alpha_n s_0}{z}u_n = 0. \qquad (3)$$

Substitution of Eqs. (2) and (3) into Eq. (1) gives

$$\sum_n a_n \left[\left(\frac{\gamma-\gamma_n}{z} + \frac{\delta}{z-1} + \varepsilon + s_0\right)u_n' + \frac{(\alpha+\alpha_n s_0)z - (q+\alpha_n s_0)}{z(z-1)}u_n\right] = 0 \qquad (4)$$

or $\quad \sum_n a_n \left[\left((\varepsilon+s_0)z(z-1) + (\gamma-\gamma_n)(z-1) + \delta z\right)u_n' + \left((\alpha+\alpha_n s_0)z - (q+\alpha_n s_0)\right)u_n\right] = 0. \qquad (5)$

To proceed further, we need recurrence relations between the involved confluent hypergeometric functions. We discuss here three different sets of such relations applying to the Kummer confluent hypergeometric functions of the form ${}_1F_1(\alpha_0+n;\gamma_0+n;s_0 z)$, ${}_1F_1(\alpha_0+n;\gamma_0;s_0 z)$ and ${}_1F_1(\alpha_0;\gamma_0+n;s_0 z)$, respectively. Consider these cases separately.

## 2. Expansions

**2.1.** $\alpha_n = \alpha_0 + n$, $\gamma_n = \gamma_0 + n$, where $n$ is an integer: $u_n = {}_1F_1(\alpha_0 + n; \gamma_0 + n; s_0 z)$.

The starting recurrence relation is simply the differentiation rule for the Kummer confluent hypergeometric functions:

$$u_n' = s_0 \frac{\alpha_n}{\gamma_n} u_{n+1}. \tag{6}$$

A further observation here is that in this case $z^2 u_n'$, $z u_n'$ and $z u_n$ are not expressed as a linear combination of functions $u_n$. If we then demand the coefficients of these terms to be equal to zero, we get that should be $s_0 = -\varepsilon$ and, furthermore, $\alpha_n = \alpha/\varepsilon$, $\gamma_n = \gamma + \delta$ for any $n$, that is the parameters $\alpha_n$, $\gamma_n$ should not depend on $n$. Thus, this contradiction leads to the conclusion that in this way no expansion can be constructed.

However, there is another possibility. Indeed, put $s_0 = -\varepsilon$ (thereby canceling the term proportional to $z^2$) and further demand

$$(\alpha + \alpha_n s_0)z - (q + \alpha_n s_0) = f(n, z),$$
$$(\gamma - \gamma_n)(z - 1) + \delta z = A f(n+1, z), \tag{7}$$

where $A$ is a constant and $f(n, z)$ is a liner function of $n$. It is then easily shown that this is possible only if $A = 1/\varepsilon$ and

$$\gamma_0 = 1 + \alpha_0 + \gamma + \delta - \alpha/\varepsilon, \tag{8}$$

$$q = \alpha - \delta\varepsilon. \tag{9}$$

Eq. (5) is then written as:

$$\sum_n a_n \left( -f(n+1, z) \frac{\alpha_n}{\gamma_n} u_{n+1} + f(n, z) u_n \right) = 0. \tag{10}$$

This gives a simple two-term recurrence relation for the coefficients of the expansion (2):

$$a_n - \frac{\alpha_{n-1}}{\gamma_{n-1}} a_{n-1} = 0, \tag{11}$$

whence $\quad u = \sum_n a_n \cdot {}_1F_1(\alpha_0 + n; \gamma_0 + n; -\varepsilon z), \quad a_n = \frac{(\alpha_0)_n}{(\gamma_0)_n}, \tag{12}$

where $(\alpha_0)_n$ is the Pochhammer symbol: $(\alpha_0)_n = \alpha_0(\alpha_0 + 1)...(\alpha_0 + n - 1)$. In general, the derived recurrence relation defines a double-sided infinite series. The series is applicable if $\varepsilon \neq 0$ and $\gamma_0$ is not zero or a negative integer. It is right-hand side terminated if $\alpha_0 = 0$ or $\alpha_0 = -N$ for some positive integer $N$. Choosing $\alpha_0 = 0$ and changing $n \to -n$ we arrive at the following expansion for $q = \alpha - \delta\varepsilon$:

$$u = \sum_{n=0}^{\infty} \frac{(1-\gamma_0)_n}{n!} {}_1F_1(-n; \gamma_0 - n; -\varepsilon z), \quad \gamma_0 = 1 + \gamma + \delta - \alpha/\varepsilon. \tag{13}$$

However, the above development can be essentially extended to avoid the additional restriction (9) imposed on the parameters of the confluent Heun equation (1). This is achieved by noting that the following recurrence relation holds

$$z(u'_n - s_0 u_n) = (\gamma_n - 1)(u_{n-1} - u_n). \tag{14}$$

Indeed, again put $s_0 = -\varepsilon$ and demand [compare with Eq. (7)]

$$\alpha + \alpha_n s_0 = -s_0[(\gamma - \gamma_n) + \delta], \tag{15}$$

which is fulfilled if

$$\alpha_0 = \gamma_0 - \gamma - \delta + \alpha/\varepsilon. \tag{16}$$

Eq. (5) is then rewritten as

$$\sum_n a_n \left[ ((\gamma - \gamma_n) + \delta) z(u'_n - s_0 u_n) - (\gamma - \gamma_n) u'_n - (q + \alpha_n s_0) u_n \right] = 0, \tag{17}$$

so that by virtue of above recurrence relations we have

$$\sum_n a_n \left[ ((\gamma - \gamma_n) + \delta)(\gamma_n - 1)(u_{n-1} - u_n) - (\gamma - \gamma_n) s_0 \frac{\alpha_n}{\gamma_n} u_{n+1} - (q + \alpha_n s_0) u_n \right] = 0. \tag{18}$$

Accordingly, the recurrence relation for the coefficients of the expansion (2) now becomes three-term:

$$R_n a_n + Q_{n-1} a_{n-1} + P_{n-2} a_{n-2} = 0, \tag{19}$$

where

$$R_n = (\gamma + \delta - \gamma_n)(\gamma_n - 1), \tag{20}$$

$$Q_n = -R_n + \varepsilon \alpha_n - q, \tag{21}$$

$$P_n = \varepsilon \frac{\alpha_n}{\gamma_n}(\gamma - \gamma_n). \tag{22}$$

For left-hand side termination of the series at $n = 0$ should be $R_0 = 0$. This is the case if

$$\gamma_0 = \gamma + \delta \quad (\Rightarrow \alpha_0 = \alpha/\varepsilon) \tag{23}$$

($\gamma_0 = 1$ is forbidden because of division by zero: $P_1 \sim 1/\gamma_{-1}$, $\gamma_{-1} = 0$). Thus, finally, the expansion is explicitly written as

$$u = \sum_{n=0}^{\infty} a_n \cdot {}_1F_1((\alpha/\varepsilon) + n; \gamma + \delta + n; -\varepsilon z) \tag{24}$$

and the coefficients of the recurrence relation (19) are explicitly given as

$$R_n = -n(\gamma + \delta + n - 1), \quad Q_n = n(\gamma + \delta + n - 1) + (\varepsilon n + \alpha) - q, \quad P_n = -\frac{(\delta + n)(\varepsilon n + \alpha)}{\gamma + \delta + n}. \quad (25)$$

Obviously, this expansion is applicable if $\varepsilon \neq 0$ and $\gamma + \delta$ is not a negative integer. A concluding remark is that the second independent solution of Eq. (1) can be constructed in a similar way after applying the transformation $z \to 1 - z$, which preserves the form of the confluent Heun equation.

The series (24) is terminated if $a_{N+1} = 0$ and $a_{N+2} = 0$ for some non-negative integer $N$. The first condition is fulfilled if $Q_N a_N + P_{N-1} a_{N-1} = 0$, while for a non-zero $a_N$ the second condition results in the equation $P_N = 0$, i.e.,

$$\alpha / \varepsilon = -N \quad \text{or} \quad \delta = -N. \quad (26)$$

The equation $a_{N+1} = 0$ (or, equivalently, $Q_N a_N + P_{N-1} a_{N-1} = 0$) then presents a polynomial equation of the degree $N+1$ for the accessory parameter $q$, defining, in general, $N+1$ values of $q$ for which the termination of the series occurs. Note, finally, that in the case $\alpha / \varepsilon = -N$ the resulting solution (24) is a polynomial in $z$.

**2.2.** $\alpha_n = \alpha_0 + n, \gamma_n = \gamma_0 = \text{const}$, where $n$ is an integer: $u_n = {}_1F_1(\alpha_0 + n; \gamma_0; s_0 z)$.

In this case we have the following recurrence relations

$$zu'_n = \alpha_n(u_{n+1} - u_n), \quad (27)$$

$$s_0 z u_n = (\alpha_n - \gamma_0) u_{n-1} + (\gamma_0 - 2\alpha_n) u_n + \alpha_n u_{n+1}. \quad (28)$$

Combining these equations we also have

$$s_0 z^2 u'_n = \alpha_n (s_0 z u_{n+1} - s_0 z u_n) =$$
$$\alpha_n \left( (\alpha_n + 1) u_{n+2} + (\gamma_0 - 3\alpha_n - 2) u_{n+1} + (3\alpha_n - 2\gamma_0 + 1) u_n + (\gamma_0 - \alpha_n) u_{n-1} \right). \quad (29)$$

Eq. (5) is rewritten as follows

$$\sum_n a_n \left[ \left( z^2(\varepsilon + s_0) + z(-\varepsilon - s_0 + \delta + \gamma - \gamma_0) + \gamma_0 - \gamma \right) u'_n + \left( (\alpha + \alpha_n s_0) z - (q + \alpha_n s_0) \right) u_n \right] = 0. \quad (30)$$

Since $u'_n$ is not expressed as a linear combination of functions $u_n$, we demand $\gamma_0 - \gamma = 0$. Now, substituting Eqs. (27), (28) and (29) into Eq. (30) we obtain a four-term recurrence relation for coefficients $a_n$:

$$R_n a_n + Q_{n-1} a_{n-1} + P_{n-2} a_{n-2} + S_{n-3} a_{n-3} = 0, \quad (31)$$

where
$$R_n = (\alpha_n - \gamma)(\alpha_n \varepsilon - \alpha), \quad (32)$$

$$Q_n = (\alpha_n \varepsilon - \alpha)(\gamma - 2\alpha_n) - s_0(\alpha_n(\varepsilon - \delta) - q) + \alpha_n(\gamma - 1 - \alpha_n)(s_0 + \varepsilon), \tag{33}$$

$$P_n = \alpha_n\left((\alpha_n + \delta)\varepsilon - \alpha + (2\alpha_n + 2 - \gamma - \delta - \varepsilon)(\varepsilon + s_0) + (\varepsilon + s_0)^2\right). \tag{34}$$

$$S_n = -\alpha_n(1 + \alpha_n)(s_0 + \varepsilon). \tag{35}$$

Note that here $s_0$ is a free parameter that can be chosen as convenient.

If we put $s_0 = -\varepsilon$ (thereby removing the $z^2$-dependence in the coefficient of $u'_n$ in Eq. (30)) then the four-term recurrence relation becomes 3-term:

$$R_n a_n + Q_{n-1} a_{n-1} + P_{n-2} a_{n-2} = 0, \tag{36}$$

where
$$R_n = (\alpha_n - \gamma)(\alpha_n - \alpha/\varepsilon), \tag{37}$$

$$Q_n = (\alpha_n - \alpha/\varepsilon)(\gamma - 2\alpha_n) + \alpha_n(\varepsilon - \delta) - q, \tag{38}$$

$$P_n = \alpha_n(\alpha_n + \delta - \alpha/\varepsilon). \tag{39}$$

The initial conditions for left-hand side termination of the derived series at $n = 0$ are $a_{-2} = a_{-1} = 0$. It then follows that should be $R_0 = 0$. This is the case only if $\alpha_0 = \alpha/\varepsilon$ or $\alpha_0 = \gamma$. Then, the final expansion is explicitly written as

$$u = \sum_{n=0}^{\infty} a_n \cdot {}_1F_1(\alpha_0 + n; \gamma; -\varepsilon z) \tag{40}$$

and the coefficients of the recurrence relation (36) take the form

$$R_n = (n + \alpha_0 - \gamma)(n + \alpha_0 - \alpha/\varepsilon), \tag{41}$$

$$Q_n = (\gamma - 2(n + \alpha_0))(n + \alpha_0 - \alpha/\varepsilon) + (n + \alpha_0)(\varepsilon - \delta) - q, \tag{42}$$

$$P_n = (n + \alpha_0)(n + \alpha_0 + \delta - \alpha/\varepsilon). \tag{43}$$

This expansion is applicable if $\varepsilon \neq 0$ and $\gamma$ is not a negative integer.

The series is right-hand side terminated at some $n = N$ if $a_N \neq 0$ and $a_{N+1} = a_{N+2} = 0$. Then, should be $P_N = 0$. If $\alpha_0 = \alpha/\varepsilon$, this condition is satisfied if

$$\alpha/\varepsilon = -N \quad \text{or} \quad \delta = -N. \tag{44}$$

If $\alpha_0 = \gamma$, the only possibility, since $\gamma$ should not be a negative integer, is

$$\gamma + \delta - \alpha/\varepsilon = -N. \tag{45}$$

Again, for each of these cases there exist $N+1$ values of $q$ for which the termination occurs. These values are determined from the equation $a_{N+1} = 0$ (or, equivalently, $Q_N a_N + P_{N-1} a_{N-1} = 0$).

**2.3.** $\alpha_n = \alpha_0$, $\gamma_n = \gamma_0 + n$, where $n$ is an integer: $u_n = {}_1F_1(\alpha_0; \gamma_0 + n; s_0 z)$.

In this case the following recurrence relations are known:

$$z u'_n = (\gamma_n - 1)(u_{n-1} - u_n), \tag{46}$$

$$u'_n = s_0 \left( u_n - \left(1 - \frac{\alpha_0}{\gamma_n}\right) u_{n+1} \right). \tag{47}$$

Since in this case $z^2 u'_n$, and $z u_n$ are not expressed as a linear combination of functions $u_n$ we demand the coefficients of these terms in Eq. (5) to be equal to zero. It is seen that should be $s_0 = -\varepsilon$ and $\alpha_0 = \alpha/\varepsilon$. So that Eq. (5) is rewritten as

$$\sum_n a_n \left[ ((\gamma - \gamma_n + \delta)z + (\gamma_n - \gamma)) u'_n + (\alpha - q) u_n \right] = 0. \tag{48}$$

Substituting Eqs. (46) and (47) into Eq. (48) we get the following recurrence relation

$$R_n a_n + Q_{n-1} a_{n-1} + P_{n-2} a_{n-2} = 0, \tag{49}$$

where
$$R_n = (\gamma + \delta - \gamma_n)(\gamma_n - 1), \tag{50}$$

$$Q_n = (\gamma + \delta - \gamma_n)(1 - \gamma_n) + \varepsilon(\gamma - \gamma_n) + \alpha - q, \tag{51}$$

$$P_n = (\gamma_n - \gamma)(\varepsilon - \alpha/\gamma_n). \tag{52}$$

For left-hand side termination of the derived series at $n = 0$ should be $R_0 = 0$. This is the case only if $\gamma_0 = \gamma + \delta$ ($\gamma_0 = 1$ is forbidden because of division by zero: $P_1 \sim 1/\gamma_{-1}$, $\gamma_{-1} = 0$). Thus, the expansion is finally written as

$$u = \sum_{n=0}^{\infty} a_n \cdot {}_1F_1(\alpha/\varepsilon; \gamma + \delta + n; -\varepsilon z) \tag{53}$$

and the coefficients of the recurrence relation (49) are explicitly given as

$$R_n = -n(\gamma + \delta + n - 1), \tag{54}$$

$$Q_n = n(\gamma + \delta + n - 1) - \varepsilon(\delta + n) + \alpha - q, \tag{55}$$

$$P_n = (\delta + n)\left(\varepsilon - \frac{\alpha}{\gamma + \delta + n}\right). \tag{56}$$

This expansion is applicable if $\alpha \neq 0$, $\varepsilon \neq 0$ and $\gamma + \delta$ is not a negative integer. If the series is right-side terminated for some non-negative integer $N$ then $P_N = 0$. This is the case if

$$\gamma + \delta - \alpha/\varepsilon = -N \quad \text{or} \quad \delta = -N. \tag{57}$$

The termination occurs for $N + 1$ values of the accessory parameter $q$ defined from the equation $a_{N+1} = 0$.

## 3. Discussion

Thus, using different recurrence relations obeyed by the Kummer confluent hypergeometric functions, we have constructed several confluent hypergeometric expansions of the solutions of the confluent Heun equation. The forms of the dependence of the used confluent hypergeometric functions on the summation variable differ from that applied in previous discussions. An example of application of the presented expansions to physical problems is the recent derivation of finite-sum closed-form solutions of the quantum two-state problem for an atom excited by a time-dependent laser field of Lorentzian shape and variable detuning providing double crossings of the frequency resonance [20].

A major set of physical problems where the presented expansions can be applied is encountered in quantum physics (see, e.g., [2-3] and references therein). For instance, in particle physics, there are many potentials for which the *stationary* Schrödinger equation is reduced to the confluent Heun equation (see, e.g., [21-22]). Similarly, in atomic, molecular and optical physics, there are many electromagnetic field configurations for which the *time-dependent* Schrödinger equations for the probability amplitudes of a driven quantum few-state system can be reduced to the confluent Heun equation [20,23-24]. Consider, for instance, the quantum two-state problem [25]:

$$i\frac{da_1}{dt} = U e^{-i\delta} a_2, \quad i\frac{da_2}{dt} = U e^{+i\delta} a_1, \qquad (58)$$

where $a_{1,2}(t)$ are the probability amplitudes of the two involved states, $U(t)$ is the Rabi frequency of the applied laser field, and $\delta(t)$ is the detuning modulation (the derivative $d\delta/dt \equiv \omega_{21} - \omega_L(t)$ is the detuning of the transition frequency from the laser frequency). This system is equivalent to the following linear second-order ordinary differential equation:

$$\ddot{a}_2 + \left(-i\dot{\delta} - \frac{\dot{U}}{U}\right)\dot{a}_2 + U^2 a_2 = 0, \qquad (59)$$

where the over-dots denote differentiation with respect to time.

Many two-level models can be solved in terms of special functions via reduction of Eq. (59) to the corresponding standard equation (see, e.g., [25]). According to the class property of solvable models, each integrable model then generates a whole class of solvable field configurations [24,26]. By applying a transformation of the dependent and independent variables one finds fifteen such four-parametric classes solvable in terms of the confluent Heun functions [20]. These classes are given as

$$U(t) = U_0 z^{k_1}(z-1)^{k_2} \frac{dz}{dt}, \quad \frac{d\delta}{dt} = \left(\delta_0 + \frac{\delta_1}{z} + \frac{\delta_2}{z-1}\right)\frac{dz}{dt}, \tag{60}$$

where $U_0$ and $\delta_{0,1,2}$ are arbitrary complex constants which should be chosen so that the functions $U(t)$ and $\delta(t)$ are real for the chosen complex-valued $z(t)$, and $k_{1,2}$ are integers or half-integers obeying the inequalities $-1 \le k_{1,2} \cup k_1 + k_2 \le 0$. The solution of the initial two-state problem is explicitly written as

$$a_2 = e^{\alpha_0 z} z^{\alpha_1}(z-1)^{\alpha_2} H_C(\gamma, \delta, \varepsilon; \alpha, q; z), \tag{61}$$

where $H_C$ is the solution of the confluent Heun equation (1), the parameters of which as well as the auxiliary parameters $\alpha_{0,1,2}$, are defined by the input parameters of the applied field.

The classes (60) extend all the previously known families of 3- and 2-parametric models solvable in terms of the hypergeometric functions to more general 4-parametric classes involving 3-parametric detuning modulation functions. Among the numerous new field configurations, several families suggest a rather convenient framework for theoretical considerations because here the parameters of the confluent Heun function are real. In several cases this allows one to derive closed form solutions using series expansions. An example of such a family is the above-mentioned generalized Lorentzian model, a family of laser pulses of Lorentzian shape with symmetric detuning functions describing both non-crossing and crossing processes with one or two level-crossing points:

$$U(t) = \frac{U_0}{1+t^2}, \quad \frac{d\delta}{dt} = \Delta_0 + \frac{\Delta_1}{1+t^2}. \tag{62}$$

The solution of the two-state problem for this model is given as

$$a_2 = z^{\alpha_1}(z-1)^{-\alpha_1} H_C(1+R, 1-R, -2\Delta_0; 0, -(R+\Delta_1/2)\Delta_0; z), \tag{63}$$

where $z = (1+it)/2$, $\alpha_1 = (\Delta_1 + 2R)/4$, and $R = \sqrt{U_0^2 + \Delta_1^2/4}$ is the effective Rabi frequency. It is readily seen that the above series (40)-(43) may terminate if $\delta = 1 - R = -N$ (i.e. if the effective Rabi frequency is a natural number: $R = 1, 2, ...$) thus resulting in closed form exact solutions. The second termination condition then defines a relation between $\Delta_0$ and $\Delta_1$ for which the termination actually occurs. Interestingly, the particular sets of the involved parameters for which these closed form solutions are obtained define curves in the 3D space of the involved physical parameters belonging to the complete return spectrum of the two-state quantum system [20]. This is readily verified using a counterpart expansion of the confluent Heun function in terms of the Tricomi confluent hypergeometric functions [27].

**Acknowledgments**

This research has been conducted within the scope of the International Associated Laboratory (CNRS-France & SCS-Armenia) IRMAS. The research has received funding from the European Union Seventh Framework Programme (FP7/2007-2013) under grant agreement No. 295025 – IPERA. The work has been supported by the Armenian State Committee of Science (SCS Grant No. 13RB-052).